\input amstex
\NoBlackBoxes



\input colordvi



\comment

olive-green  7,0,95,35       green (51,0,96,39)       red (0,100,87,14)
    darkblue (94,100,0,32)   blue (98,30,0,28)      darkbrown (0,32,61,82)
    brown (0,43,83,58)      purplered (5,98,0,29)         purple (44,97,0,44)
     verydarkblue(96,97,0,84)     pink (0,28,6,4)         orange (0,53,99,31)
\endcomment




\font\rm=cmr10 \rm

\font\bf=cmb10
\font\Rm=cmr9 at 11pt
\rm
\font\it=cmsl9 at 10pt
 at 7pt

\font\Rrm=cmr17 at 16pt
   \font\Rm=cmr12 at 11.5pt

\long\def\Pf{\par\noindent {\it Proof.} }
\def\({\left(}
\def\){\right)}
\def\st{such that }
\def\qed{\hfill$\bullet$\vskip 4pt}

\def\brcs#1{\left\{ #1\right\}}

\def\Log{\text{Log\,}}
\def\iso{\cong}
\def\wrt{with respect to }
\def\:{\,:}

\def\supp{\text{supp}\,}

\def\T{\text{\bf T}}

\def\Re{{\text{Re}\,}}

\def\Exp{\text{Exp\,}}

\def\R{\text{\bf R}}
\def\N{\text{\bf N}}
\def\Z{\text{\bf Z}}

\def\Arrow #1;#2.{#1\:#2 \to }

\def\Set#1#2{\brcs{#1 \left|\vphantom{#1 #2} \right.#2}}

\def\Oh#1{{\pmb O}\(#1\)}

\def\oh#1{{\pmb o}\(#1\)}

 \def\supp#1{\text{supp\,}#1}
\def\Rrr#1,#2{{\Cal J}_{#1,#2}}
\def\slfrac#1#2{{\raise -.07 ex\hbox{$^{#1}$}}\!/\raise .35 ex \hbox{${}_{#2}$}}
\def\ssf #1/#2{\slfrac {#1}{#2}}

\def\pd #1,#2.{\frac {\partial #1}{\partial #2}}

   \long\def\Lem
#1.#2\par{\vskip4pt{\baselineskip=13pt\font\it=cmsl12 at
11.5pt\Rm
   \noindent {\rm \uppercase{#1}} #2\vskip3pt

   }} 

\long\def\Proclaim #1.#2 \endproclaim{\vskip4pt{\baselineskip=13pt\font\it=cmsl12 at
11.5pt\Rm
   \noindent {\rm \uppercase{#1}} #2\vskip3pt

   }} 

\long\def\remark #1\endremark{\vskip 2pt \noindent {\it Remark\/} #1\par}

\long\def\Sectionhead #1.#2:\par #3{\vskip 4pt \noindent {\bf #1 #2}vskip 2pt\noindent\nospace #3}

\long\def\Title #1\par {\noindent{\Rrm #1}\vskip 9pt}

 \long\def\SubT #1.{\noindent {\it #1\/} } 
 
 \long\def\SecT
#1\par{\vskip 3pt \noindent {\bf #1}\vglue1pt
   \noindent}

\long\def\subtitle #1.{\vskip 2pt \noindent {\it #1}}

\long\def\Rmk#1\par{\vskip 1pt \noindent {\it
Remark.} #1\vskip2pt}

\long\def\Abstract #1\par{{\leftskip= 3 true cm \rightskip = 3 true cm \font\it=cmsl10 \font\rm=cmr10 \baselineskip = 10pt
\parindent=.35 true cm\rm\noindent 
{\it Abstract} #1\vskip 8pt

}}

\long\def\Author #1 \par{\noindent{\it #1}}
 
\input diagrams.tex

\def\onezer{0.0}
\def\oneone{0.1}
\def\onetwo{0.2}
\def\onethr{1.1}
\def\onefou{1.2}
\def\onefiv{1.3}
\def\onesix{1.4}
\def\onesev{1.5}
\def\oneeig{1.6}
\def\onenin{1.7}
\def\oneten{1.8}
\def\oneele{1.9}
\def\onetwe{1.10}
\def\onethi{1.11}
\def\oneftn{1.12}

\def\twoone{2.1}
\def\twotwo{2.2}
\def\twothr{2.3}
\def\twofou{2.4}
\def\twofiv{2.5}

\let\iso=\cong

\def\paren#1{\/{\rm(#1)}}

\def\tripnorm #1xxx{\left\|\hglue-.2ex\left|#1\right|\hglue-.2ex\right\|}

\def\Exp{\text{Exp}\,}

\def\Log{\text{Log\,}}

\def\op{{}^{\text{op}}}

\def\SS{{\Cal S}}
\def\MM{{\Cal M}}
\def\NN{{\Cal N}}

 \Title Another invariant for AT actions %

\Abstract  
 We construct a collection of numerical invariants for approximately transitive (AT) actions (of $\Z$). We use them (sometimes supplemented by other invariants to show that members of various one-parameter families of AT actions are mutually non-isomorphic. 

{}\plainfootnote{}{AMS(MOS) classification: 28D05, 37A99, 46L10;%
 key words and phrases:  approximately transitive, tensor product of actions, variance
}

\Author David Handelman 

\SecT Introduction

Let $(X,\mu)$ be a measure space, and let $\Arrow T; X.X$ measurable invertible ergodic transformation.The classification of $(T,X,\mu)$ \wrt measurable conjugacy was shown ([CW]) to be equivalent to the classification (up to isomorphism) of their   von Neumann algebra crossed products, $L^{\infty}(X) \times_T \Z$, and in turn of matrix-valued random walks in terms of a boundary (called the Poisson boundary). This was subsequently ([GH]) shown to be equivalent to classification of a measure-theoretic version of dimension groups. 

Approximately transitive (AT) actions are those that can be expressed in either of the latter two formulations as the $1 \times 1$ matrix case or, respectively, rank one (not in the ergodic sense, which is far too strong). In particular, AT actions correspond to  direct limits, subject to an equivalence relation (that corresponds to conjugacy, isomorphism, etc, in the other formulations). We now describe this.

Let $(P_m)$ be a sequence of members of $l^1 (\Z)$ (which we usually view as functions on the unit circle) \st $P$ has only nonnegative coefficients and  $P_m(1) = 1$. We may form the direct limit
$$\diagram
l^1(\Z)  & \rTo^{\times P_1 }  &   l^1(\Z) &\rTo^{\times P_2} & l^1(\Z)  & \rTo^{\times P_3} &  l^1(\Z)  & \rTo^{\times P_4} & \dots   ,
\enddiagram
$$
where $\Arrow \times P_m; l^1(\Z) . l^1(\Z)$, sending the $m$th copy of $l^1(\Z)$ to the $m+1$st, is given by multiplication of functions, that is, $f \mapsto P_m \cdot f$ (equivalently, by convolution with the distribution corresponding to $P_i$). We also assume that the infinite product  of the $P_m$ with any translation by powers of $x$, does {\it not\/} exist (this is to guarantee  nontrivial ergodicity in this context). There is a lot of structure preserved by the maps (such as positivity), and then there is a completion process. 

Fortunately, we do not have to go over this, because isomorphism is described by a relatively simple equivalence relation involving almost commuting diagrams.  Let $(P_m)$ and $(Q_m)$ be sequences as above. The following diagram describes the equivalence relation.  There exist  telescopings, $P^{(u(i))} :=  P_{u(i)}P_{u(i) + 1}\cdots P_{u(i+1)-1} $,  and $Q^{v(i)} := Q_{v(i)}Q_{v(i) + 1}\cdots Q_{v(i+1)-1} $ (the functions $u$ and $v$ demarcate the telescopings).

$$
\diagram
 \dots & l^1(\Z) &   \rTo^{\times P^{(u(i))} }  &   l^1(\Z) &  \rTo^{\times P^{(u(i+1))} }  & l^1(\Z)  & \rTo^{\times P^{(u(i+2))}} &  l^1(\Z)  & \rTo^{\times P^{(u(i+3))}} & \dots  \\
&\dTo^{R_i} &   \!\!\!\!\!S_i{{\nearrow}} & \dTo^{R_{i+1}}& {\!\!\!\!\!S_{i+1}} {\nearrow}& \dTo^{R_{i+2}}  &  {\!\!\!\!\!S_{i+2}}{\nearrow}&\dTo^{R_{i+3}}  & {\!\!S_{i+3}}{\nearrow}&\dTo^{} & \\
 \dots & l^1(\Z) &   \rTo^{\times Q^{(v(i))} }  &   l^1(\Z) &  \rTo^{\times Q^{(v(i+1))} }  & l^1(\Z)  & \rTo^{\times Q^{(v(i+2))}} &  l^1(\Z)  & \rTo^{\times Q^{(v(i+3))}} & \dots  \\
\enddiagram
$$
Here $R_i, S_j$ are elements of $l^1(X)$ with only nonnegative coefficients, and $R_i (1) = S_j (1) = 1$ (we can reduce to the case that additionally, $R_i$ and $S_j$ are Laurent polynomials---the corresponding distributions are finite), and multiplication by $R_i$ sends the $u(i)$th copy of $l^1(\Z)$ in the top row to the $v(i)$th copy of $l^1(\Z)$ in the bottom  row, while $S_j$ sends the $v(j)$ copy in the bottom row to the $u(j+1)$st copy in the top row. The functions $u$ and $v$ serve as index functions. 

In particular, $R_{j+1} S_j$ sends the $v(j)$th copy on the bottom row to the $v(j+1)$st (on the same row). But we have another obvious map that does this, specifically, the product $Q_{v(j+1)-1} Q_{v(j+1)-2}\cdots Q_{v(j)}$, which we have denoted $Q^{(v(j))}$. Then we require the summability condition
$$
\sum_j \left\| R_{j+1}S_j - Q^{(j)}\right\| < \infty.
$$
(The $l^1 (\Z)$-norm is used.)

Similarly, $S_i R_i$ sends the $u(i)$th copy of $l^1(\Z)$ on the top row to the $u(i+1)$. Set $P^{(i)}$ to be the corresponding product, $P_{u(i+1)-1}P_{u(i+1)-2}\dots P_{u(i)+1}P_{u(i)}$. Then we require that 
$$
\sum_i \left\| R_{i}S_i - P^{(u(i))}\right\| < \infty.
$$

The existence of  $R_i, S_j$ satisfying all these conditions is equivalent to there being an isomorphism between the von Neumann algebras, or conjugacy of the corresponding ergodic transformations, etc [GH, Theorem\, 3.1]. (So we don't have to know how to define the Poisson boundary, for example.)

If $\MM$ is the (isomorphism class of)  von Neumann algebra corresponding to the sequence $(P_m)$ (or any sequence equivalent to it), we typically write, {\it $\MM$ corresponds to the system $(P_m)$\/} or vice versa. 

While this yields the correct notion of isomorphism, it is rarely easy to decide on isomorphism or nonisomorphism of two systems using it. Invariants have been developped. The best known and earliest is the T-set, corresponding to eigenvalues of the transformation on a suitable Banach space. This is relatively easy to calculate, but only coarsely separates systems (algebras). A massive family of numerical invariants was introduced in [GH] and used there and in [H], which we will call {\it mass-cancellation invariants\/} (they will be defined and used in section 2). An unpublished result of Giordano, Handelman, and Munteanu asserts that the T-set invariants can be recovered from the mass-cancellation invariants. 

Mass cancellation invariants are often useful, but very often, are difficult to calculate. In this paper, we introduce a family of numerical invariants that are generally easier to calculate. We use them to show that for many natural one-parameter families of AT actions (more precisely, their von Neumann algebras), $(\MM(r))_{r \in \R^{++}}$, the members, $\MM(r)$ are mutually non-isomorphic. 

On the other hand, mass cancellation invariants can distinguish (in many cases) an AT transformation from its inverse---which the new invariants cannot. 

Section 0 describes the new invariant, and presents some elementary properties. Section 1 contains applications to divisible systems (where the $P_m$ are compound Poisson), culminating in the nonisomorphism Theorem \oneftn. Section 2 contains applications to not necessarily divisible systems. Here, there are more complications, necessitating that mass cancellation invariants  assist in distinguishing systems.

Let $(P_m)_{m \in \N}$ be a sequence of Laurent polynomials or absolutely summable Laurent series in one variable, with only nonnegative coefficients, and such that $P_m (1) = 1$ for all $m$. This describes an AT (approximately transitive) action, although there is no guarantee in this generality that it is nontrivial. 

 We describe a numerical invariant for equivalence (measure-theoretic isomorphism) that often allows to distinguish members of one-parameter families of these. 
 For example, let $r$ be a positive real number and  let $n$ be a positive integer exceeding $1$. Let $\MM (r)$  denote the AT action arising from $\(\Exp (rx^{n^m})\)$; here $P_m (x) = \exp \(rx^{n^m} -1 \)$. It is fairly well known  that $\MM (r) \iso \MM (r')$ entails $N = N'$. Using the new invariant, we provide more results of this type. 

\SecT 0 The invariant

 It is not one invariant, but an uncountable collection of numerical invariants, analogous to those in [GH, H]. Let $(w_k)_{k \in \N}$ be a sequence of complex numbers of modulus $1$. There is no restriction on the sequence. We wish to associate a number, denoted $\SS ((w_k), (P_m))$, in $[0,1]$, so that the assignment $(P_m) \mapsto \SS ((w_j), (P_m))$ is an isomorphism invariant for $(P_m)$. In other words, each sequence $(w_k)$ yields an isomorphism invariant for AT actions. Most sequences of elements of the circle yield uninteresting or simply uncomputable invariants, but for the examples we have in mind, there are natural choices of sequences which yield nonisomorphism results. 

 Fix the sequences $(w_k)$ and $(P_m)$.  In direct analogy  with the invariant discussed in [GH,  H], define for each  $l \in \N$, the number $S_{k,l}$ defined as 
$$
 S_{k,l} = \lim_{d \to \infty}  \left|\prod_{m=l}^{m= l+d}P_m (w_k) \right| . \tag1
 $$
 That the limit exists follows from $|P(z)| \leq 1$ when $P$ has no negative coefficients, $P(1) = 1$, and $|z| = 1$. Now define 
$$
 S_l = \inf_{k \in \N} S_{k,l}.
$$
 Finally, set $ \SS ((w_j), (P_k)) = \lim_l S_{l}$. That the limit exists follows from $S_l \leq S_{l+1}$. This is a number in the unit interval, but there is no guarantee that it is nonzero or not $1$. 

Our first task is to show that this is indeed an isomorphism invariant. This is routine, but presented here to convince skeptics (and better than Arthur Cayley's proof of the Cayley-Hamilton theorem.)

\Lem Proposition \onezer. Suppose that $\MM_1$ corresponds to the system $(P_m)$ and $\MM_2$ corresponds to $(Q_m)$. Let $(w_k)$ be a sequence of elements of the unit circle. If $\MM_1 \iso \MM_2$, then $\SS((w_k), (P_m)) = \SS((w_k), (Q_m))$. 

\Pf We are given the almost commuting diagram given above. Given $\epsilon > 0$, there exists $j'$  \st for all $j \geq j$, 
$$\eqalign{
\left\|  R_j S_{j-1} R_{j-1}S_{j-2} \cdots R_{j' + 1}S_{j'} - Q^{(v(j))}Q^{(v(j-1))}\cdots Q^{(v(j'))}\right\| & < \epsilon \cr
\left\|   S_{j-1} R_{j-1}S_{j-2} \cdots S_{j'+1}R_{j' + 1} - P^{(u(j))}P^{(u(j-1))}\cdots Q^{(u(j'))}\right\| & < \epsilon \cr
}$$
Now $|R_j S_j\cdots S_{j'}(w_k)| \leq |S_j \cdots R_{j'+1} (w_k)|$ (the second product is simply the first with the first and last terms deleted). As $| Q^{(v(j))}Q^{(v(j-1))}\cdots Q^{(v(j'))}(w_k) -  R_j S_{j-1} R_{j-1}S_{j-2} \cdots R_{j' + 1}S_{j'}(w_k) | < \epsilon$ and similarly with the truncated version, we obtain $S_{k, u(j')} \leq S^Q_{k,v(j')} + \epsilon$ for all sufficiently large $j'$. It follows that $S_{u(j')}  \leq S^Q_{v(j')} + \epsilon$ for infinitely many $j'$, and thus $\SS((w_k), (P_m)) \leq \SS((w_k), (Q_m)) +\epsilon$. As this is true for all $\epsilon$, we have $\SS((w_k), (P_m)) \leq \SS((w_k), (Q_m))$. Reversing the roles of $P_m$ and $Q_m$, we obtain the opposite inequality, so $\SS((w_k), (P_m)) = \SS((w_k), (Q_m))$.
\qed

So we can write $\SS((w_k), \MM)$ for $\SS((w_k), (P_m))$ if $(P_m)$ corresponds to $\MM$. If we replace $\inf$ by $\sup$ in the definition, we obtain another invariant ({\it upper\/}; so the first invariant is the {\it lower\/} one associated to $(w_k)$; but in most of the examples here, the upper and lower ones agree. We never use the  upper invariant.

 As noted before, it is similar to the mass cancellation invariants introduced in [GH]. It has one advantage over these, in that it is usually easier to compute with, at least if we make the appropriate choice for the sequence $(w_k)$. 

 \noindent{\it Tensor products of actions.} If $\MM$ and $\NN$ are the von Neumann algebra crossed products associated to  two actions, then we may form their von Neumann algebra tensor product, $\MM \otimes\NN$. If the actions are AT, say corresponding to $(P_m)$  and $(Q_m)$ respectively, then the action arising from the sequence of products $(P_m Q_m)$  corresponds to $\MM \otimes \NN$. However, it is not clear what the dynamical interpretation should be; that is, if $T$ and $U$ are ergodic transformations, how should $T \otimes U$ be defined? A first guess is $T \times U$, but this need not be ergodic. An attempt to resolve this,  intended for the topological, rather than measure-theoretic, setting (minimal replacing ergodic) is given in [BH, Appendix A].

The invariant is frequently (but not always)  multiplicative  \wrt tensor products, that is, if $\MM$ and $\MM'$ are AT actions, then $\SS((w_j), \MM \otimes \MM') = \SS((w_j), \MM ) \cdot  \SS((w_j),  \MM' ) $ often occurs.

\Lem Lemma \oneone. Let $\MM$ and $\MM'$ be AT, and let $(w_k)$ be a sequence of elements of the unit circle. 
\item{(a)} $\min \brcs{\SS ((w_k), \MM)  , \SS ((w_k),\MM') } \geq \SS ((w_k), \MM \otimes \MM') \geq \SS ((w_k), \MM)  \cdot \SS ((w_k),\MM') $;
\item{(b)} If $\SS ((w_k), \MM) \neq 0$, then $\SS ((w_k)), \MM \otimes \MM') = 0$ iff $\SS ((w_k),   \MM') = 0$;
\item{(c)} $\SS ((w_k), \MM \otimes \MM)  = \SS ((w_k), \MM)^2$;
\item{(d)}  If $\SS ((w_k), \MM) =1$, then $\SS ((w_k), \MM \otimes \MM') =\SS ((w_k)),\MM') $.

\Pf If $\MM$ is given by the sequence $(P_m)$ and $\MM'$ is given by $(P'_m)$, then $\MM \otimes \MM'$ is given by $(P_m P'_m)$ (the sequence of products). If $S_{k,l}$, $S'_{k,l}$, and $S''_{k,l}$ represent the moduli of the products for $(P_m)$, $(P'_m)$ and $(P_m P'_m)$ respectively, then we have $S''_{k,l} = S_{k,l}\cdot S'_{k,l}$. Taking infima over $k$, we obtain $S''_l \geq S_{l}\cdot S'_{l}$. Taking limits as $l \to \infty$, the right inequality of (a) follows. 

Next, we see that $S''_{k,l} =  S_{k,l}\cdot S'_{k,l} \leq \min \brcs{S_{k,l},S'_{k,l}}$, and the left side of (a) follows. 

\noindent (b) If $\SS ((w_k)),   \MM') = 0$, then $\SS ((w_k)), \MM \otimes \MM') = 0$ follows from the left side of (a). If $\SS ((w_k)),   \MM') \neq 0$, then $\SS ((w_k)), \MM \otimes \MM') \neq 0$ follows from the right side of (a). 

\noindent (c) Follows from $|P_m (w_k)^2| = |P_m(w_k)|^2$.

\noindent (d) Follows from both parts of (a). 
\qed

The following construction is   obvious, but is given here for completeness. 

\Lem Example \onetwo. Two AT systems, $\MM_1$ and $\MM_2$, and a sequence of roots of unity, $(w_k)$, \st $\SS((w_k), \MM_1 \otimes \MM_2) \neq \SS((w_k), \MM_1)\cdot \SS((w_k), \MM_2)$.

\Pf Let $\tau$ be a positive real number less than $\slfrac12$, and let $w_k = \exp (2 i \pi/2^k)$ Define 
$$
P_m = \cases \frac{1 +2x^{2^m}}3 & \text {if $m$ is even} \\
\frac{1 +\tau x^{2^m}}{1+\tau} & \text {if $m$ is odd.}
\endcases
$$
Let $\MM_1$ denote the AT action determined by $(P_m)$. If we interchange {\it even\/} with {\it odd} in the definition of $P_m$, we obtain $\MM_2$. It is routine to verify that, for the infinite products, 
$$\eqalign{
\(1 - \frac 89\)&\(1- \frac{4\tau}{(1+\tau)^2}\cdot \frac 12\) \(1 -\frac 89 \sin^2 \frac{2\pi}{16} \) \(1 -\frac{4\tau}{(1+\tau)^2}  \sin^2 \frac{2\pi}{32}\) \cdot  \dots \cr 
& <\(1- \frac{4\tau}{(1+\tau)^2}\) \(1 -\frac 89 \cdot \frac12\) \(1 -\frac{4\tau}{(1+\tau)^2}  \sin^2 \frac{2\pi}{16}\) \cdot  \dots .\cr
}$$

This implies that for every $l$, $S_l $ (for $\MM_1$) is the square root of the top product, hence this is $\SS((w_k)), \MM)$. The same inequality yields that this is also $\SS((w_k), \MM_2)$. 

However, $\MM_1 \otimes \MM_2$ is given by the sequence $(Q_m = (1+ 2x^{2^m})(1 + \tau x^{2^m})/3(1+\tau))$, and it is easy to check that 
$$
\SS((w_k)), \MM_1 \otimes \MM_2)   = \(\prod \( 1- \frac 89 \sin^2 \frac{2\pi}{2^m}\) \cdot \prod \( 1- \frac {4\tau}{(1+\tau)^2} \sin^2 \frac{2\pi}{2^m}\)\)^{1/2},
$$
which is not $\SS((w_k)), \MM_1) \cdot \SS((w_k)), \MM_2) = \SS((w_k)), \MM_1)^2$.
\qed

 If $\MM \otimes \MM \iso \MM$ (as occurs for many odometers), then the values of the new invariants can only be $0$ or $1$, no matter what the choice of sequence $(w_k)$. A little more generally, if $\MM(r)$ (for $r \in \R^{++}$) is a one-parameter family of AT actions \st $\MM(r) \otimes \MM(r') \iso \MM(r+r')$ and the invariants $\SS ( (w_k),\cdot)$ are multiplicative on $\brcs{\MM(r)}$ (the original example, above, of $\MM (r)$, satisfies these properties), then for fixed $(w_k)$, the map $\phi:\R^{++} \to [0,1]$ given by $r \mapsto \SS ((w_k), \MM(r))$ satisfies $\phi(r+r') = \phi(r)\phi(r')$ (e.g., $\phi(r) = \gamma^r$ for some $\gamma \leq 1$). In this case, if $\phi(r) \not\in \brcs{0,1}$ for some $r$, then $\MM(r) \iso \MM(r')$ entails $r  = r'$ (Corollary \onefou), which is precisely the conclusion we want.

\SecT 1 Bounded AT actions

Here we give reasonably general sufficient conditions so that one-parameter families of AT actions, $\MM(r)$ (for $r$ a positive real number), given by, for example,  $(P_{m,r} = \exp\(r(h_m (x^{n^m}-1) \)$ (for some positive integer $n \geq 2$) where each $h_m \in l^1(\Z)$ has no negative coefficients and $h_m(1) = 1$, satisfy $\MM(r') \iso \MM(r)$ implies $r = r'$. The corresponding distributions are compound Poisson, and  therefore divisible.

\comment

The criterion ($\dag$) in the following is satisfied generically,  and yields lots of examples for which the invariant is multiplicative. The proof is elementary. 

Do we need both satisfy dagger?

\Lem Lemma W R O N G. Suppose that $\MM_1$ and $\MM_2$ are AT systems, and $(w_k)$ is a sequence of elements of $\T$. Suppose that   $\MM_1$ is given by $(P_m)$, and the corresponding $S_{k,l}$ satisfy {\par}
\noindent $(\dag)$ for infinitely many $l$, $S_l = \liminf_{k} S_{k,l}$. {\par}\noindent Then $\SS((w_k), \MM_1 \otimes \MM_2) = \SS((w_k), \MM_1)\cdot \SS((w_k),  \MM_2)$. 

\Pf Let $\MM_2$ be given by $(Q_m)$, so that the tensor product is given by $(P_m Q_m)$. Pick on of the $l$s for which ($\dag$) is satisfied. Given $\epsilon > 0$, there exist infinitely many $k$ \st  $S_{k,l} \leq S_l + \epsilon$. Let $S_{k,l}'$ be the counterparts for $(Q_m)$, and $S_{k,l}''$ for $(P_mQ_m)$. From 
$S_{k,l}'' = S_{k,l}S_{k,l}'$, we have for infinitely many $ k$, $S_{k,l}'' \leq (S_l + \epsilon) S_{k,l}'$ 

\endcomment

 Let $f = \sum_{t \in \Z}  a_t x^t$ be an element of $l^1(\Z)$, that is, $\sum |a_t| < \infty$, \st all the coefficients, $a_t$, are nonnegative. In that case, $\| f \| = f(1)$. In general, for $f \in l^1(\Z)$), $\| f \| \geq \sup_{z \in \T} |f(z)|$. We say $f$ has {\it finite second moment,} if $f$ has only nonnegative coefficients, $f(1) = 1$, and $\sum a_t t^2  < \infty$. This implies that  the real and imaginary parts of $f$ (as a function on the unit circle) are twice differentiable. When $f$ has finite second moment, we define $\mu_2 (f):= \sum a_t t^2$;   this is $f' (1) + f'' (1)$. The {\it first moment\/}, $\mu_1(f) =\sum  a_t t = f'(1)$, is defined if merely $\sum a_t |t| < \infty$, but need not be positive. 
 
 First we discuss the most general divisible AT situation.  Let $H_m$ be elements of $l^1(\Z)$ with no negative coefficients. Define $P_m= \Exp H_m  := \exp (H_m - H_m(1))$, and let $\MM$ denote the system corresponding to the sequence $(P_m)$. We first note that we can assume each $H_m$ is a Laurent polynomial (that is, has finite support), by a perturbation result (that is, the system so obtained is isomorphic to $\MM$).  In this generality, there is no guarantee that the action is nontrivial.

  For real  $r > 0$, define $\MM (r)$ to be the action obtained from the sequence $(P_{m,r} := \Exp (r H_m)$ ($= \exp (r (H_m - H_m(1)$). The following is elementary. 

 \Lem Proposition \onethr. Let $\MM$ be the action  determined by $(\Exp H_m)$. Suppose that $\SS ((w_k), \MM) = S$   for some choice of sequence of elements of $\T$, $(w_k)$. Then for $r >0$, we have $\SS ((w_k), \MM(r)) = S^r$. 

 \Rmk Our convention is that $S^r$  is $\exp (r \ln S)$ if $S > 0$, and $0$ if $S = 0$. 
 
 \Pf If $w$  is of absolute value $1$, then 
$$\eqalign{
 \left| P_{m,r} (w) \right| & = \exp \(\Re (r(H_{m,r} (w) -H(1)) ) \) \cr
 & = \(\exp (\Re(H_m(w)-H(1))\)^r.\cr 
 }$$
It follows that $S_{k,l, r} := \lim_{d \to \infty} \prod_{j = l}^{l+d} |P_{j,r} (w_k)|$ is the exponential (in $r$) of   $S_{k,l,1}$. Thus for fixed $r$, the infimum over $k$ of $S_{k,l,r}$ is just the exponential of the corresponding number with $r=1$ (because $r \mapsto \gamma^r$ is order preserving). Now taking limits (as $l \to \infty$), the result follows. \qed 

Recall that for $H \in l^1({\Z})$ with only nonnegative coefficients and $H(1)  = 1$,  $\Exp H$ denotes $\exp (H-1) \in l^1(\Z)$.

\Lem  Corollary \onefou. Let $(H_m)$ be a sequence of elements of $l^1(\Z)$ with only nonnegative coefficients, and \st $H_m(1) = 1$ for all $m$. Let $r > 0$, and let  $\MM(r)$ be the AT action determined by $(\Exp (r (H_m))$. If for some choice of sequence of elements of $\T$, $(w_k)$, and some $r' > 0$, we have $\SS((w_k),\MM(1)) \neq 0,1$, then ${\MM(r)}$ are mutually non-isomorphic. 

\Pf Without loss of generality, we may assume $r' = 1$. Let $S = \SS((w_k),\MM(1)) $; then by the preceding, $\SS((w_k),\MM(r)) = S^r$, and of course $r \mapsto S^r$ is one to one. \qed

Thus, to determine whether $\MM(r) \iso \MM(r')$ implies $r = r'$, it is sufficient to show that $S:= \SS((w_k), (\MM))$ is neither zero or one. This is   useful, as computing $S$ exactly can be quite difficult.

 Let $f \in \l^1(\Z)$ be \st $f(1) = 1$ and $f$ is C$^2$, that is, $\mu_2 (f) <\infty$. Define $V(f)$ to be $f''(1) + f'(1) - (f'(1))^2$. It is straightforward that for $f_1$ and $f_2$, we have $V(f_1 f_2) = V(f_1) + V(f_2)$, and if $f = x^n$ for some integer $n$, then $V(f) = 0$. In particular, $V(x^n f) = V(f)$, that is, $V$ is shift invariant. 

 Of course, $V$ is well-known. Let $(a_j)_{i\in \Z}$ be a sequence of nonnegative real numbers \st $\sum a_j= 1$, and $\sum a_j j^2 < \infty$, and let $f = \sum a_j x^j$; then $V(f)$ is simply the {\it variance\/} of the distribution $(a_j)$, or equivalently of the corresponding integer-valued random variable $X$ defined by $\text{Pr} (X= j) = a_j$.   In our definition of $V$, there is no requirement that the $a_j$ be nonnegative or even real.  We can similarly define (unnormalized) skewness and kurtosis (provided the third and fourth moments respectively are finite), as well as   third and fourth cumulants; the latter convert multiplication (of functions) to sums (as variance does), but it is unlikely they will be of any use here. 

 The following three elementary results (Lemmas \onefiv, \onesix, and Corollary \onesev) are probably known, but I could not find   references for them. The variance results will be useful in section 2. 
 
 \Lem Lemma \onefiv. Let $f = \sum a_j x^j$ have finite second moment, with $f(1) = 1$. Then 
$$V(f) = \sum_{j < j'} a_j a_{j'} (j'-j)^2. 
$$ 

 \Pf Everything is absolutely summable here, so there will be no problem with the infinite sums. We  observe 
 $$\eqalign{ \mu_2 (f) = 1 \cdot \mu_2 (f) &= \(\sum_j a_j\)\cdot \(\sum_t a_t t^2\) \cr
 & = \sum a_j^2 + \sum_{j \neq j'} a_j a_j' (j')^2  \cr
 & = \sum a_j^2 + \frac 12 \sum_{j \neq j'} a_j a_{j'}( j^2 + (j')^2  )\cr
& =  \sum a_j^2 + \sum_{j < j'} a_ja_{j'} (j^2 + (j')^2). \cr 
 V(f) & = f''(1) + f'(1) -( f'(1))^2\cr
 & = \mu_2 (f) - \(\sum a_j j  \)^2 \cr
 & = \mu_2 (f) - \sum a_j^2 j^2  - 2 \sum_{j<j'} a_j a_{j'}  jj'. \quad \text{Thus}\cr
 \sum_{j < j'} a_j a_{j'} (j'-j)^2 & =  \sum_{j < j'} a_j a_{j'} (j^2 + (j')^2)  - 2 \sum_{j < j'} a_j a_{j'}jj' \cr 
 & = \( \mu_2 (f) - \sum a_j^2 j^2\) + V(f) - \( \mu_2 (f) - \sum a_j^2 j^2\)\cr
 & = V(f). \cr
}$$ \qed

The following generalizes the fact that for $\theta > 0$, we have $\sin \theta < \theta$. The proof of the following is by double induction, observing that the derivative transforms the partial sums into shorter versions of their ilk. 

\Lem Lemma \onesix. Suppose $\theta > 0$. Then 
\item{(i)} for all nonnegative integers $s$, $$
 \sum_{t=0}^{2s+1}  (-1)^t \frac{\theta^{2t+1}}{(2t+1)!}<\sin \theta < \sum_{t=0}^{2s} (-1)^t \frac{\theta^{2t+1}}{(2t+1)!};
$$
\item{(ii)} for all $s\geq 1$, $$
 \sum_{t=0}^{2s+1} (-1)^t \frac{\theta^{2t}}{(2t)!} <\cos \theta < \sum_{t=0}^{2s} (-1)^t \frac{\theta^{2t}}{(2t)!}.$$ 

In the following, $K(h) = \sum_{j < j'} a_j a_{j'} (j'-j)^4/12$. Unlike the case with $4$  replaced by $2$ (which yields a multiple of variance), this does not seem to be very interesting (and we never use it). 

\Lem Corollary \onesev. Let $h \in l^1 (\Z)$ have only nonnegative coefficients, and  $h(1) = 1$; assume the fourth moment exists. Let $\theta> 0$ and set $z = e^{i\theta}$. 
\item{(i)} $-\frac{\mu_2 (h)}{2}\theta^2 + \frac{\mu_4 (h)}{24}\theta^4  > \Re (h(z) -1) > -\frac{\mu_2 (h)}{2}\theta^2$
\item{(ii)} $1 - V(h) \theta^2 + K(h) \theta^4 >|h(z)|^2 > 1 - V(h)\theta^2$.

\Pf  Write $h = \sum a_j x^j$ with $a_j \geq 0$, $\sum a_j = 1$, and $\sum a_j j^4 < \infty$. Then 
$$\eqalign{ 
\Re (h(z) -1) &= \sum a_j \cos j\theta -\sum a_j \cr
& = \sum  a_j (\cos j\theta -1);\quad \text{by Lemma \onesix,}\cr
-\sum  a_j \(\frac{(j\theta)^2 }2 - \frac{(j\theta)^4}{24} \) &  > \Re (h(z) -1) > -\sum a_j \frac{(j\theta)^2 }2;  \quad \text{so}\cr
-\theta^2\frac{\mu_2 (h)}2 + \theta^4 \frac{\mu_4 (h)}{24} & > \Re (h(z) -1) > -\theta^2\frac{\mu_2 (h)}2.\cr
}$$

We also have 
$$\eqalign{
|h(z)| ^2& = \left| \sum a_j (\cos j\theta + i \sin \theta)\right|^2\cr
& = \sum a_j^2 + 2\sum_{j < j'} a_j a_{j'}\((\cos j\theta) \cdot (\cos j'\theta) + (\sin j\theta) \cdot (\sin j'\theta)\)\cr
&= \sum a_j^2 + 2\sum_{j < j'} a_j a_{j'}(\cos (j'-j)\theta); \quad \text{since $1 = \sum a_j^2 + 2\sum_{j< j'} a_j a_{j'}$,}\cr
&= 1 -  2\sum_{j < j'} a_j a_{j'}(1- \cos (j'-j)\theta); \quad \text{thus by  Lemma \onesix,}\cr
}$$
$$
\eqalign{
1- \sum_{j< j'} a_j a_{j'}(j'-j)^2 \theta^2 +2 \sum_{j< j'} a_j a_{j'}\frac{( (j'-j) \theta)^4 }{24} &>  |h(z)| ^2  > 1 -  \sum_{j< j'} a_j a_{j'}(j'-j)^2 \theta^2;  \quad  \text{so}\cr
1- V(h) \theta^2 + K(f)\theta^4 &>  |h(z)| ^2  > 1- V(h) \theta^2. \cr
}$$
\qed

 If in this lemma, $f$ is a Laurent polynomial, say with $M$ as maximal exponent and $m$ as minimal one, we can use the Bhatia-Davis inequality [BD], $V(f) \leq (M - f'(1))(f'(1)- m)$, to bound $V(f)$   without going to the trouble of  calculating it.

 Let ${(n(m))}_{m = 1,2, \dots}$ be a sequence of positive integers, and form $T(m) = \prod_{l=1}^{m} n(l)$. In the simplest case, $n(m) = n$ for all $m$, and then $T(m) = n^m$. Now let $(h_m)$ be a sequence of elements of elements of $l^1 (\Z)$, each with no negative coefficients and zero constant term satisfying $h_m (1) = 1$, and each of finite second moment. Let $r$ be a positive real number, and set $P_{m,r} (x) = \Exp \(rh_m(x^{T(m)}\)$; this is $\exp \(r(h_m (x^{T^m} -1)\)$, owing to the normalization of the $h_m$. We   form the AT action, $\MM (r)$, given by the sequence $P_{m,r}$, yielding a one-parameter family, $r \mapsto \MM(r)$. The current aim is to determine sufficient conditions so that $\MM(r) \iso \MM(r')$ implies $r = r'$, and for this, it is sufficient to show $\SS((w_k),\MM) \not\in \brcs{0,1}$, by Corollary \onefou. 
 
 For example, this property fails   if for some $r$, we have that $\MM(r)$ is an odometer, for then $\MM(2r)  \iso \MM(r) \otimes \MM (r) \iso \MM(r)$. The simplest example of this occurs if $n(m) = n > 1$ for all $m$ and $h_m (x)= f(m) x $ where  $f$ is a positive valued function \st $f(m)/\ln m \to \infty$. Then  $P_{m,r} = \Exp (rf(m) x^{n^m})$, and by [H, Theorem 4.4], $\MM (r)$ is the $n$-odometer for all choices of $r > 0$; the corresponding supernatural number is $n^{\infty}$, that is, infinite at the prime divisors of $n$ and zero at other primes.
 
  If we weaken the hypothesis on the growth of $f$, it is known [op\.cit.] that sufficient for the system to be an odometer is that $\liminf f(m)/\ln m$ be sufficiently large (how large depending on $n$), but it is not known whether $\MM(r)$ is an odometer if $f$ grows more slowly, e.g., $f(m) \sim \ln\ln m$. Unfortunately, our new invariant doesn't help with this. However, it will help if, for example, $h_m = h$ for all $m$, and some variations on this, e.g., a bound on second moments of $h_m$.  

We can now deal with the relatively simple case of $h_k = h$ for all $k$. First, the case of $T(k) = n^k$, that is, $n(k) = n$ is constant. 

\Lem Proposition \oneeig. Let $n \geq 2$ be an integer, and  $h  \in l^{1}(\Z)$ have only nonnegative coefficients,  finite second moment, and $h(1) = 1$. Set  $\MM (r)$ to be the system associated to $\(P_m := \Exp (rh(x^{n^m}))\)$. Let $w_k = \exp (2 \pi i/n^k)$. Then 
$$ 
S:= \SS((w(k), \MM(1))) =\exp\(- \sum_{t=1}^{\infty } \Re (1- h(e^{2\pi i/n^t}))\),
$$
and this lies strictly between $0$ and $1$. In particular, $\brcs{\MM(r)}_{r > 0}$ are mutually nonisomorphic. 

\Rmk We can get fairly tight estimates for  $S$ if the first few terms of the sum are known. 

\Pf That $S$ is as given in the display is an immediate consequence of the definitions. To check that $S$ is neither zero nor one, we note that 
$$
\Re (1- h(e^{2\pi i/n^t})) \leq  2 \mu_2(h)\frac {\pi^2}{n^{2t}},
$$
hence the sum converges, and thus $S > 0$; but $S < e^{-\Re (1- h(\exp 2\pi i)/n)}$ (the first term), and so $S < 1$. The rest follows from Corollary \onefou. 
\qed

The remark follows from the inequality $1 - \mu_2 (h) \theta^2 /2< \Re (1- h(e^{i\theta }) <  1 - \mu_2 (h) \theta^2/2 + \mu_4(h) \theta^4/24$, which for $\theta = 2\pi/n^t$ gives a   tiny error  for sufficiently large $t$. 

If $n(k) \to \infty$ (but $h_k = h$ for all $k$), then the invariant does nothing; the value will be $1$. This will follow from more general results,  where   $h_k$ are allowed to vary. 
 
Returning to our general situation (with $n(m)$, $T(m)$, $h_m$ being general), let $w_k =\exp (2\pi i/T(k))$, a primitive $T(k)$th root of unity. We will determine (under relatively modest conditions) the value of $S:= \SS ((w_k), (P_{m,1}) )$. We are particularly interested in sufficient conditions so that it is neither zero nor one, as $\SS((w_k), (P_{m,r}) = S^r$ (this means the positive value of $s^r$, $\exp(r\ln S)$, when $s \neq 0$).  

 Abbreviate $P_{m,1}$ to $P_m$. First, we note that $|P_{m,r} (z)| = |P_m (z)|^r$ for $z$ on the unit circle. Next, we see that $P_m(w_k) = 1$ if $k \leq m$. For $k > m$, set $\theta_{m,k} = \exp (2\pi i/(T(k)/T(m)))$.  We have 
$$\eqalign{
 |P_m (w_k)| & = \left|\exp \( h\(\exp \frac{2\pi i}{T(k)/T(m)}\) -1 \) \right| \cr
 & = \exp \(-\Re (1 -h (\exp (i\theta_{m,k})) \).\cr 
}$$ 
 Thus, for $l $, $d$ positive integers, we have 
$$\eqalign{
  \prod_{j=l}^{l+d}  \left| P_{j} (w_k)\right| & = \prod_{j=l }^{(l+d) \wedge (k -1)}  \left| P_{j} (w_k)\right|; \quad \text{thus, as $d \to \infty$},\cr 
 & = \exp \( -\sum_{j = l}^{k-1} \Re (1 - h_j(e^{i\theta_{j,k}} ))\); \quad \text{substituting $t = k-j$,} \cr 
 & =  \exp \( -\sum_{t =1}^{k-l} \Re (1 - h_{k-t}(e^{i\theta_{j,k} }))\) \cr 
 & = \exp (- \Re (1-h_{k-1} (\exp (2\pi i/n(k))) - \Re (1- h_{k-2} (\exp (2\pi i/n(k)n(k-1))) - \dots )) \cr 
 }\tag*$$ 
 The last line is purely expository.

The following yields conditions under which the value of the invariant is not zero.

 \Lem Lemma \onenin. Suppose that $\mu_2 (h_k) < \infty$ for all $k$, and in addition, 
 \item{(a)}  $$\limsup \frac{\mu_2 (h_{k-1})}{\mu_2(h_k)) n(k)^2}:= C < 1;$$ 
 \item{(b)} $$\limsup \frac{\mu_2(h_{k-1})}{n(k)^2 n(k-1)^2} := \rho < \infty.$$ 
 {

} \noindent Then, with $ \MM$ given by $\(P_m (x) = \Exp (h_m (x^{T(m)}))\)$, we have 
$$
 \SS ((w_k), \MM) \geq \exp \(- \limsup_{k\to\infty} \Re (1 - h_k(e^{2\pi i/n(k+1)}))\) \exp (- M\rho),
$$
where     $M = 2\pi^2/(1-C)$, with equality if $\rho = 0$. 
 
 \Rmk In particular, the value of the invariant is not zero here.  

\Rmk Hypothesis (a) is quite weak. Hypothesis (b), $\mu_2 (h_{k-1}) = \Oh{n(k)^2 n(k-1)^2}$, is reasonable.

 \Pf  Expression  (*) yields a value for the product, and thus for $S_{k,l}$. For the  purposes of simplicity of the terms, set $n(0) = 1$. In the penultimate line thereof, take the sum beginning with $t = 2$; define 
\comment 
 Applying the mean value theorem, there exists $\theta* \in [0,\theta]$ \st $H(\theta) - H(0) = \theta H'(\theta*)$. Thus
$$\eqalign{
 H'(\theta*) &= \sum a_j j  (1 - \cos j\theta*)\cr 
  &\leq \sum a_j( j^3 \theta*^2/2 )  
}$$
 Since $H(0) = 1$, we thus have $1- H(\theta) \leq \mu_3 (h) \theta^3$ (or $d^2 \mu_1(h) \theta^3/2$ where $d$ is the degree of $h$). 
 \endcomment
 
$$\eqalign{
 A_k&:= \Re \(1-h_{k-2} (e^{2\pi i/n(k)n(k-1)})\) +\cdots + \Re\(1-h_l(e^{2\pi i/n(k)\cdots n(l+1)}) \)\cr
 &\leq \sum_{t=2} ^{k-l}  \(\frac{2\pi}{n(k) n(k-1)\cdots n(k-t+1)} \)^2  \mu_2 (h_{k-t})/2. \cr
& = \frac{2\pi^2}{n(k)^2n(k-1)^2} \sum_{t=2} \frac{\mu_2(h_{k-t}) }{(n(k-2)  \cdots    n(k-t+1))^2}\cr
 }$$  

 Now let $C = \limsup \mu_2 (h_{k-1})/(\mu_2 (h_k)\cdot n(k)^2)$, and pick $C'$ \st $1 >C' > C$. There exists $k_0$ \st $j \geq k_0$ entails $\mu_2 (h_{j-1}) \leq C' \mu_2 (h_j) n(j)^2$. Iterating this when   $l > k_0$, we obtain
$$\eqalign{
 A_k &\leq \frac{2\pi^2}{n(k)^2n(k-1)^2} \sum_{t=2} \mu_2(h_{k-1}) (C') ^{t-2}
\cr
& \leq 2\pi^2\frac{\mu_2(h_{k-1})}{n(k)^2n(k-1)^2} \frac{1}{1-C'}. \cr
}$$

Now $S_{k,l} = \exp(-\Re(1- h_k(e^{2\pi i/n(k+1)})) - A_k$, hence for sufficiently large $  l$,
$$
S_{k,l} \geq \exp(-\Re(1- h_k(e^{2\pi i/n(k+1)})) \cdot e^{-M\rho}
$$
where $M = 2\pi^2/(1-C')$.
 
If $\rho = 0$, we obtain $S_{l}  \geq \exp \(- \limsup \Re (1 - h_k(e^{2\pi i/n(k+1)}))\)$ (for sufficiently large $l$), and the reverse inequality is trivial. \qed  

Now we can (almost) finish the $h_m = h$ case. 

\Lem Corollary \oneten.   Let $n(k) \to \infty$  and  $h  \in l^{1}(\Z)$ have only nonnegative coefficients,  finite second moment, and $h(1) = 1$. Set  $\MM (r)$ to be the system associated to $\(P_m := \Exp (rh(x^{T(m)})\)$. Let $w_k = \exp (2 \pi i/n(k+1))$. Then 
$$ 
\SS((w(k), \MM(r))) =1.
$$

\Pf Without loss of generality, we can assume $r =1$. Conditions (i) and (ii) are satisfied with $C = \rho = 0$, yielding $\SS((w(k), \MM(r))) = \exp \(- \limsup \Re (1- h(e^{2\pi i/n(k+1)})\)$, but this is clearly $1$. \qed

\comment
Residual cases for $h = h_k$: (a) $ n(k)$ bounded; the same argument that works for constant $n(k)$ works here, and (b) $\sup n(k) = \infty$ but there is an infinite subset of $\N$, $Y$, \st $\brcs{ n(k)}_{k \in Y}$ is bounded. This probably behaves as in case (a), but I have not verified it. 

\endcomment

 Now we obtain estimates for $\Re (1 -h_{k-1}(e^{2\pi i/n(k)})$; it  is equivalent, and slightly more convenient to work with $\Re (1 -h_{k}(e^{2\pi i/n(k+1)})$. 
 
 So let $h = \sum a_j x^j$ with $a_j\geq 0$ for all $j$, and $\sum a_j = 1$. Let $n$ be a positive integer exceeding $1$, and let $R >1$ be a real number. Let $\supp h$ denote the set of $j$ \st $a_j \neq 0$.

 Define 
$$\eqalign{
 S(h,n,R) &= \supp h \cap \(\bigcup_{t \in \Z}  \( tn + \left[\frac nR , n\cdot (1 - 1/R)\right]\)\) \cr 
 U(h,n,R) & = \sum_{j \in S(h,n,r)} a_j
}$$

For example, if $\theta = \pi/n$, then $j \in tn + [n/R, n(1-1/R) $ for some integer $t$  entails that $j\theta \in t\pi + [\pi/R, \pi (1-1/R)]$, and thus $\sin^2 j\theta \geq \sin^2 \pi/R$. On an interval of the form $[0, K]$, the proportion of it not in the union,  is about $2n/R$, so for large $R$,  $S(h,n,R)$ is typically most of $\supp h$. If the distribution of $h$ is not concentrated off $S(h,n,R)$, then $U(h,n,R)$ will be close to $1$, or at any rate, more than one-half. If we can arrange that this occurs uniformly in $k$ for $h_k $ and $n(k+1)$ (playing the roles of $h, n$ respectively) for some $R$, then we obtain a lower bound for values of the invariants.

 \Lem Lemma \oneele.  Suppose there exists $R > 1$ \st $\liminf_k U(h_k, n(k+1), R) := \eta > 0$. Then for all sufficiently large $k$, $\Re (1 - h_k(e^{2\pi i/n(k+1)})) \geq 2 \eta \sin^2 (\pi/R)$. 

 \Pf Write $h_k = \sum_j a_{j,k}$, so that 
$$
 \Re(1 - h_k(e^{i\theta}) = \sum_j a_{j,k} (1- \cos j\theta) = 2 \sum_j a_{j,k}\sin^2 (j\theta/2).
$$ 
 Set  $\theta = 2\pi/n(k+1)$. We see that for $j \in S(h_k, n(k+1), R)$, we have $\sin^2 (j\theta/2) \geq \sin^2 \pi/R$. Hence for all sufficiently large $k$, $\Re (1 - h_k(e^{2\pi i/n(k+1)})) \geq 2 \eta \sin^2 (\pi/R)$. \qed

 There follows immediately: 
 
 \Lem Corollary \onetwe. Suppose that there exists $R > 1$ \st for all sufficiently large $k$, there exists $\eta > 0$  \st $U(h_k, n(k+1),R) \geq \eta $. Then $\SS ((w_k), \MM) \leq e^{-2\eta \sin^2 (\pi/R)}$.  
 
 In particular, this yields a fairly weak  sufficient condition (on the sequence $(h_m)$) so that the value of the invariant is strictly less than $1$. 

We also have a converse to this. 

 \Lem Proposition \onethi. Suppose  that     for all $R > 1$, $\liminf_k U(h_k, n(k+1), R) = 0$.  Then $\SS ((w_k), \MM) = 1$. 

 \Pf  Let $h = \sum  a_j x^j$, and $\theta$ a small positive real number.   For $j \not\in S(h,n, R)$, $|\sin j\theta| < \pi/R$. Thus
$$\eqalign{
 \sum_{j \not\in S(h,n,R)} a_j (1-\cos 2j\theta) 
& = 2\sum_{j \not\in S(h,n,R)} a_j \sin^2 (j\theta); \cr
& <\(\frac {\pi}R\)^2. 
}$$
Therefore
$$
\Re (1- h(e^{2i\theta}) = \Re \sum a_j (1-\cos 2j \theta) < U(h,n,R) + \(\frac {\pi}R\)^2, 
$$
and so sufficient for the   left side to be small is that both summands be small. 

Suppose we have $U(h_k,n(k+1),R) < \epsilon$ for infinitely many $k$, and let $\theta = \pi/n(k+1)$. Since $\mu_2(h) = \oh{n(k+1)^2}$, $\Re (1-h_k (e^{2\pi i/n(k+1)})$, for infinitely many $k$, $\Re(1-h_k (e^{2\pi/n(k+1)})) < \epsilon + (\pi/R)^2$. Allowing $R\to \infty$ and  $\epsilon \to 0$, we deduce $S_l \to 1$ along infinitely many $l$, and thus $\SS ((w_k), \MM) = 1$. \qed   
 
\Lem Theorem \oneftn. Let $\MM(r)$ be given by $\(P_{m,r} = \Exp (rh_m (x^{T(m)})\)$, subject to the following conditions. 
\item{(a)}$\mu_2 (h_k) < \infty$  for all but finitely many $k$;
\item{(b)} $\limsup \frac{\mu_2 (h_{k-1})}{\mu_2(h_k)) n(k)^2}  < 1$;
\item{(c)} $\mu_2 (h_{k-1}) = \Oh{n(k)^2 n(k-1)^2}$. 
\item{(d)} There exists $R > 1$ \st $\liminf_k U(h_k, n(k+1), R)  > 0$.
{\par \noindent Then  $\MM(r) \iso \MM(r') $ implies $r = r'$. 

}

\Pf Let $\MM$ denote $\MM(1)$. Set $w_k = \exp(2\pi i/T(k)$. By Lemma \onenin, $\SS((w_k),\MM) > 0$, and by Corollary \onetwe, $\SS((w_k),\MM) < 1 $. Corollary \onefou\ allows to conclude. 
\qed

Hypothesis (a) obviously holds if  the $h_m$ are Laurent polynomials; (b) is a very weak condition; and (c) is somewhat restrictive (and implies (a)), but it is difficult to see how it could be weakened. Hypothesis (d) is not very strong, but is superficially complicated.  

\comment
\Lem Proposition. Suppose  that $\mu_2 (h_k) = \oh{n(k+1)^2}$ and   for all $R > 1$, $\liminf_k U(h_k, n(k+1, R)) = 0$.  Then $\SS ((w_k), \MM) = 1$. 

 \Pf  Now consider, for $h = \sum  a_j x^j$, 
$$\eqalign{
 \sum_{j \not\in S(h,n,R)} a_j (1-\cos j\theta) 
& \leq \sum_{j \not\in S(h,n,R)} a_j j^2 \theta^2\cr
& < \mu_2 (h) \theta^2.
}$$
Thus 
$$
\Re (1- h(e^{i\theta}) = \Re \sum a_j (1-\cos j \theta) \leq U(h,n,R) +\mu_2(h) \theta^2, 
$$
and so sufficient for the   left side to be small is that both summands are small. 
Suppose we have $U(h_k,n(k+1),R) < \epsilon$ for infinitely many $k$, and let $\theta = 2\pi/n(k+1)$. Since $\mu_2(h) = \oh{n(k+1)^2}$, $\Re (1-h_k (e^{2\pi i/n(k+1)})$, for infinitely many $k$, $\Re(1-h_k (e^{2\pi/n(k+1)})) \leq \epsilon + \oh{1}_k$. Since this is true for all $\epsilon > 0$, we deduce $S_l \to 1$ along infinitely many $l$, and thus $\SS ((w_k), \MM) = 1$. \qed   
\endcomment

 \SecT 2 A different type of one-parameter family 

 In the cases discussed earlier, the mapping (for appropriate choice of $(w_k)$) $r \mapsto \SS((w_k), \MM(r))$ is multiplicative, that is, $\SS ((w_k), \MM(r) \otimes \MM(r'))) = \SS((w_k), \MM(r)) \cdot \SS((w_k), \MM(r'))$. There is another, fairly natural type of one-parameter family, for which similar properties do not apply, but nonetheless, we can obtain similar isomorphism results. 

 Let $h_m \in l^1(\Z)$ with only nonnegative coefficients, and this time we drop the normalization constraint. For each positive real $r$, define   $P_{m,r} (x) = h_{m} (rx^{T(m)})/h_m(r)$. The system now need not be divisible (as it was in the earlier case, owing to the definition of $\Exp$). It {\it can} be divisible; for example, with $n(k) = n$ and $h_m = 1 + x + \cdots + x^{n-1}$,  we see that $(P_{m,1})$ is just the $n$-odometer. 

 To distinguish this construction from the earlier ones, we use the notation $\NN$ (or $\NN(r)$) for the system arising from $(h_{m} (rx^{T(m)})/h_m(r))$.

A more natural definition might seem to be $P'_{m,r} = h_m ((rx)^{T (m)})$, but this often results in atoms, for example, if $h_m = (1+x)/2$ and $r \neq 1$, for the resulting sequence $(P'_{m,r})$, the products actually converge, resulting in  an atomic dynamical system.

\Lem Lemma \twoone. Let $h \in l^1(\Z)$ have only nonnegative coefficients, finite second moment, and $h(1) = 1$. Let $n\geq 2$ be a positive integer, and $r $ a positive real number. Let $\NN$ be of the form $(P_m)$, where $P_m = h(rx^{n^m})/h(r)$. Set $w_k = \exp 2\pi/n^k$. Then for all $l$,  $\SS((w_k),\NN) = \lim_{k\to \infty} S_{k,l}$, and this equals $\prod_{t=1}^{\infty} \(| h(r \exp (2\pi i/n^t))|/h(r)\)$. Moreover, this is nonzero unless for some $t$, $h(r\exp (2\pi i/n^t)) = 0$.  

\Pf Let $h_r (x) = h(rx)/h(r)$. First, from Lemma \onenin(ii), we have that $  |h(re^{i\theta}|^2 \geq  1- V(h_r) \theta^2$. With $\theta $ equalling successively $2\pi/n^m$, we see that $1 \geq |h(re^{2\pi i/n(k+1)\cdots n(l+1))}|^2/h(r)^2 \geq 1 - V(h_r)\pi^2/(n(k+1)\cdots n(l+1))^2$. Thus  $\prod_{t=1}^{\infty}  \( | h(r \exp (2\pi i/n^t))|/h(r)\)$ converges  in the sense of infinite products, and the only way the limit can be zero is if one of the factors is. 

We have $S_{k,l} = \prod_{t=1}^{k-l}\( |h(re^{2\pi i/n^{t}})| /h(r)\)$; fixing $l$ and taking the infimum over $k$, noting that $k-l \to \infty$, we simply obtain $S_l = \prod_{t=1}^{\infty}  \(| h(r \exp (2\pi i/n^t))|/h(r)\)$. As this is independent of $l$, we obtain $\SS((w_k),\NN) = S_l$. 
\qed

 It is not true that $\NN(r + r') \iso \NN(r) \otimes \NN(r')$ (except under  degenerate circumstances), so that multiplicativity is not as  interesting as in the previous class of examples. 

 Asking the same question, can the  class of evaluation invariants distinguish members of $\brcs{\NN(r)}$, the answer is somewhat different---it often can, but with the occasional aid of another invariant.  The following simple-looking example illustrates  what can happen. 
 
\Lem Example \twotwo. A  one-parameter family $\NN(r)$ \st $\SS((w_k), \cdot)$ distinguishes $\NN(r)$ from $\NN(r')$ if $r' \neq {r,r^{-1}}$. An additional invariant distinguishes $\NN(r)$ from $\NN(r^{-1})$ if $r \neq 1$.

 \Pf Set $h_m = 1 + x$, so that $P_{m,r} = (1+ rx^{T(m)})/(1+ r)$. If $n (k) = 2$ for all $k$, when $r=1$, the corresponding system is the dyadic odometer---but for all other values of $r$, it isn't an odometer (the former statement is elementary, the latter is not difficult, and will follow from the computation of the invariant anyway). 

 Taking our usual $w_k = \exp (2\pi i/T(k))$, we compute enough of the invariant to obtain a slightly limited classification result, which will later be supplemented by another invariant. 

 An elementary computation reveals that for $k > m$,
 $$\left| P_{m,r}(w_k)\right| ^2 = 1 -  \frac{4r}{(1+r)^2} \sin^2 \frac{\pi}{T(k)/T(m)}. 
 $$ Thus 
$$
 \left| \prod_{t = 0}^{k+l-1} P_{l+t,r} (w_k) \right|^2 = \prod_{t = 0}^{k+l-1} \( 1 -  \frac{4r}{(1+r)^2} \sin^2 \frac{\pi}{T(k)/T(l+t)}\) .
 $$ 
 The smallest term in this product is $1 -  \sin^2 (\pi/n(l+1))4r/(1+r)^2$, and it is easy to check that the product converges (as we let $k \to \infty$) in the usual sense of infinite products---however,  some of the initial terms might turn out to be zero (this occurs with the odometer example), so that the product could be zero. The product is invariant under $r \mapsto r^{-1}$ and here, $\NN(r^{-1})$ corresponds to the inverse transformation to $\NN(r)$. Thus  $\SS((w_k), \NN(r)) = \SS((w_k), \NN(r^{-1})) $. In particular, the invariant does not distinguish some pairs of members of the family. We will deal with this shortly. 

 We will show that $r \mapsto \SS((w_k), \NN(r))$ is monotone decreasing on $(0,1]$, and strictly decreasing under mild assumptions. The latter entails members of this part of the family are mutually nonisomorphic. Then we will show that if $r \neq 1$, $\NN(r)\not\iso\NN(r^{-1})$ by an easy application of the invariants introduced in [GH]. 
 
 The function $r \mapsto  4rc/(1+r)^2$  ($c$  a positive constant) is strictly  increasing on $(0,1]$ with  maximum at $r = 1$. A minor problem arises when a few factors in the  product, 
$$
 \alpha(k,r) := \prod_{t=1}^{\infty} \(1 - \frac{4r}{ (1+r)^2} \sin^2 \(\frac{\pi}{n(k+1)n(k)\cdots n(l-+2-t)}\)\),
$$
 might be zero. First, we observe that each term is nonnegative (since $4r/(1+r)^2 \leq 1$). Thus the value  zero can only occur if $r =1$ and $n(k+1) = 2$. Since $n(k) = 2$ for all $k$ and $r = 1$ entails the system corresponds to the dyadic odometer, we can set this case aside. In particular, $r \neq 1$ entails each term is positive. Moreover, since $r \mapsto 4r/(1+r)^2$ is strictly increasing, we see that for $r < r' < 1$, we have $S_{k,l} (r) > S_{k.l} (r')$. A   consequence is that $\SS((w_k), \NN(r)) \geq \SS((w_k), \NN(r')) $, but we want strict inequality. 

 This does not always hold (as we will see, when we discuss the condition $n(k) \to \infty$). However, if $n(k) = n$ (for all but finitely many $n$, then we easily see that 
$$
\SS((w_k), \NN(r)) =\prod_{j=1}^{\infty} \( 1 - \frac{4r}{(1+r)^2}\sin^2 \(\pi/n^j\)\).
$$
 The infinite product converges to a nonzero positive number, and it is easy to see that it is strictly increasing as $r \to 1$  from below. 

 Hence if $n(k) = n$, and for $r, r' \in \R^{++}$ with $r,r^{-1} \neq r'$, then $\NN(r) \not\iso \NN(r')$. 

 Now we show that $\NN(r) \not \iso \NN(r^{-1})$ if $r \neq 1$. We use the other class of invariants, $\S((p_k), (P_n))$, defined ([GH]) as follows. The $p_k$ are elements of $l^1 (\Z)$ of norm one, and we define for each $k \geq l$, $s_{l,k} = \lim_{d \to \infty}\| p_k P_l \cdot P_{l+1} \cdots P_{l+d} \|$ (the limit exists since the sequence is monotone decreasing; all the $P_k$ have norm $1$). Then define $s_l = \inf_k s_{l.k}$, and as with the evaluation invariant, note that the sequence $(s_l)$ is increasing, and thus $\lim_l s_l$ exists. This is the invariant associated to the sequence $(p_k)$, and to distinguish this type from the other one, we refer to the former as   {\it mass-loss invariants.} 

 We sometimes abbreviate $P_{k,r} = (1+r x^{T(k)})/(1+r)$ to $P_k$, if $r$ is understood.  

\Lem Lemma \twothr. Let $\NN(r) $ be given by $(P_{k,r})$. If either of the following hold, then  $\NN(r) \not \iso \NN(r^{-1})$ when $r\neq 1$.
 \item{(i)}$n(k) > 2$ for all but finitely many $k$
\item{(ii)} $n(k) = 2$ for all but finitely many $k$. 

\Rmk The mixed case, that $\Set{k \in\N}{n(k) =2}$ is both infinite and co-infinite in $\N$, is problematic. 

\Pf  If $h$ is in $l^1(\Z)$, we denote by $h\op$ the element of $l^1$ given by $x \mapsto x^{-1}$, that is, all exponents are replaced by their negatives. When $h$ is a polynomial (that is, has support in $\Z^+$ and this is finite), we can replace $h\op$ by $x^d h\op$, where $d$ is the degree of $h$, and so continue to work with polynomials (rather than Laurent polynomials).  Since $h \mapsto h\op$ is an isometry of $l^1(\Z)$ preserving all the coefficients (just reflecting them), we see immediately that 
$$
\S((p_k)), (P_m\op)) =\S((p_k\op)), (P_m)).
$$ 

Hence to show that $\S((p_k)), (P_m\op)) \neq \S((p_k)), (P_m))$ for a choice of $p_k$, it  is sufficient to show that   $\S((p_k)), (P_m)) \neq \S((p_k\op)), (P_m))$.

\vskip2pt
\noindent (i) {\it $n(k)\geq 3$ for all but finitely many $k$.} Set $p_k = (1-rx^{T(k)})/(1+r)$. We notice that $p_k \cdot P_k = (1-r^2x^{2T(k)})/(1+r)^2$, which has norm $(1-r^2)/(1+r)^2 = (1-r)/(1+r)$. Assuming $k > l$, multiply this by $P_l \cdots P_{k-1}$. This has total degree $T(l) + T(l+1) + \dots T(k-1) < T(k)$, so that the largest difference between exponents is less that $T(k)$. It follows that there is no further mass cancellation in the product $P_l \cdot P_{k-1}(p_k \cdot P_k )$, that is, $\| p_k P_l \cdots P_k\| = (1-r)/(1+r)$. 
 
 Moreover, the product polynomial $p_k P_l \cdots P_k$ has degree $T(l) + T(l+1) + \dots + T(k-1) + 2 T(k) < T(k+1)$ (this uses $n(k ) \geq 3$). Any product of the form $P_{k+1} \cdot \dots \cdot P_{l+d}$ (with $d > k-l$) is supported on $T(k+1)\Z$, so that 
 $$\eqalign{
 s_{k.l} &=   \lim_{d \to \infty}\| p_k P_l \cdot P_{l+1} \cdots P_{l+d} \| \cr 
 & = \lim_{d \to \infty}\| (p_k P_l \cdot P_{l+1} \cdots P_k) (P_{k+1} \cdots P_{l+d})\| \cr 
 & =  \lim_{d \to \infty}\| (p_k P_l \cdot P_{l+1} \cdots P_k\| \cr 
 & = \frac{1-r}{1+r}
}$$
 Thus $\S((p_k), (P_m)) = (1-r)/(1+r)$.

 By the preliminary comment, we are   reduced to showing $\S((q_k), (P_m) ) \neq \S((p_k), (P_m) $ where $q_k = (x^{T(k)}-r)/(1+r)$. Now $q_k P_k = (x^{2T(k)} +(1-r) x^{T(k)} - r)/(1+r)^2$. This has norm $2(1-r)/(1+r)^2$. Now let $Q$ be any polynomial with only nonnegative coefficients. In order for $\| q_k P_k Q\| < \| Q \| \cdot \| q_k P_k\|$, there must exist two points in the support of $Q$ whose difference is either $T(k)$ or $2 T(k)$. But no such exists in a polynomial $Q = P_l \cdots P_{k+1} \cdot P_{k+1} \cdot P_{l+d}$. Hence, as in the previous case, $\S((q_k), (P_m)) = 2(1-r)/(1+r)^2$. This is not equal to $(1-r)/(1+r)$, unless $r=1$.  
 
\noindent (ii) {\it $n(k)= 2$ for all but finitely many $k$.}
Special techniques are needed to deal with non-noninteractivity. We require some preliminary results.

For a nonnegative integer $j$, let $\delta (j)$ denote the number of $1$s in its binary expansion, and if $j \neq 0$, let $e(j)$ be the maximum power of $2$ that divides $j$. Thus $e(j ) = 0$ iff $j$ is odd, $e(j) = 1$ iff $j \equiv 2 \pmod 4$, and so on. 

Let $r$ be a positive real number, and form the product of polynomials in the variable $X$,
$$\eqalign{
Q(X) &:= \prod_{i=0}^{d-1} (1 + rX^{2^i}) \cr
& =  \sum_{j=0}^{2^d-1} r^{\delta(j)}X^j. \cr
}$$
The last line follows easily from uniqueness of   binary expansions. Evaluating at $X =1$, we obtain $\sum r^{\delta(j)} = (1+r)^d $. Let $a$ be a positive real number, and consider the product, 
$$
(1- aX)\cdot Q = 1 + \sum_{j=1} ^{2^d -1}\(r^{\delta(j)} - a r^{\delta(j-1)}\)X^j - a r^{d}X^{2^d}.
$$
In order the compute the $l^1$-norm of this, we observe that for $1 \leq j < 2^d$, we have $\delta(j-1) = \delta(j) -1 + e(j)$. 

\Lem Lemma \twofou.  For $u = 0,1,\dots, d$, the following holds. 
$$
\sum_{\Set{1 \leq j \leq 2^d-1}{e(j)= u}} r^{\delta(j)} = r (1+r)^{d-1-u}.
$$

\Rmk Of course, this is consistent with the earlier expansion, since $r\sum_{u=0}^{d-1 } (1+r)^u = (1+r)^d -1$.

\Pf Fix $u$; then $e(j )= u$ means $j \equiv 2^{u} \pmod {2^{u+1}}$ (even when $u = 0$). For $u = 0$, we have $\delta(j) = \delta(j-1)+1$, so the sum on the left becomes $r\sum r^{\delta(j)}$ where now $j$ varies over all the even integers less than or equal to $2^d-1 $. But the latter sum is the same as the sum over all terms up to $2^{d-1}$, hence is just $(1+r)^{d-1}$. 

For $ u \geq 1$, all $j$ are divisible by $2^u$, and the quotient is odd. Thus applying the result of the previous paragraph (to odd integers less than $2^{d-1-u}$), we obtain the result. 
\qed

\noindent{\it Resumption of proof of Lemma 2.3\paren{ii}.} We have 
$$\eqalign{
(1- aX)\cdot Q &= 1 + \sum_{j=1} ^{2^d -1}\(r^{\delta(j)} - a r^{\delta(j-1)}\)X^j - a r^{d}X^{2^d}\cr
&= 1 - a r^{d}X^{2^d}+  \sum_{j=1} ^{2^d -1} r^{\delta(j)} X^j(1- a r^{e(j)-1})\cr 
& = 1 - a r^{d}X^{2^d} + \sum_{u=0}^{d-1}\sum_{\Set{1 \leq j \geq 2^d-1}{e(j)= u}} X^j r^{\delta (j)} (1- a r^{e(j)-1}); \quad\text{so}\cr 
\|(1- aX)\cdot Q \|&= 1 + ar^d + \sum_{u=0}^{d-1} \sum_{\Set{1 \leq j \geq 2^d-1}{e(j)= u}} r^{\delta (j)} |1- a r^{e(j)-1}|\cr
& = 1 + ar^d + \sum_{u=0}^{d-1} r (1+r)^{d-u-1} |1- a r^{e(j)-1}|.\cr
 }$$\ 

In the special case that $r < 1$ and $a =r$, then $1-ar^{u-1} =1- r^{u}$; this is nonnegative, so that 
$$\eqalign{
\|(1- rX)\cdot Q \|& = 1  + r\sum_{u=0}^{d-1}  (1+r)^{d-u-1}(1- r^u) + r^{d+1}\cr
& = 1 + r\sum_{u=0}^{d-1}  (1+r)^{d-u-1} - r\sum_{u=0}^{d-1}  (1+r)^{d-u-1}r^u + r^{d+1}\cr
& = (1+r)^d - r (1+r)^{d-1} \sum_{u=0}^{d-1}\( \frac{r}{1+r}\)^u + r^{d+1}\cr
&= (1+r)^d - r (1+r)^{d-1} \frac{1 - \(\frac r{r+1}\)^d}{1- \frac r{r+1}}  + r^{d+1}\cr
& = (1+r)^d (1 - r) \( 1- \(\frac r{r+1}\)^d\) + r^{d+1}; \quad \text{normalizing,}\cr
\left\|\frac{1- rX)}{1+r} \cdot \frac{Q}{(1+r)^d} \right\| & = \frac{1-r}{1+r} +\Oh{ \(\frac{r}{1+r}\)^{d}} \cr
}$$

Now let $a= 1/r$, so that $1 - ar^{u-1} = 1  - r^{u-2}$. This is negative  for 
 for $u= 0,1$, zero for $u = 2$, and positive for $u \geq 3$. In that case, we obtain 
$$\eqalign{
\left\|(1- r^{-1}X)\cdot Q \right\|& = 1 + r(1+r)^{d-1} (r^{-2} -1) + r(1+r)^{d-2} (r^{-1}-1) + r\sum_{u\geq 3} (1+r)^{d-u-1}(1-r^{u-2}) + r^{d-1} \cr
& = 1 + r^{d-1}+ r(1+r)^d\cdot \cr 
& \qquad\times\( r^{-2}(1-r )+ r^{-1}\frac{1-r}{(1+r)^2} + \frac{1}{1+r}\sum_{u \geq 3} \frac{1}{(1+r)^u} - \frac{1}{(1+r)r^2}\sum_{u \geq 3}\(\frac r{1+r}\)^u \)   \cr
& =1 + r^{d-1}+  r(1+r)^d \cdot \cr 
&\qquad \times \( r^{-2}(1-r )+ r^{-1}\frac{1-r}{(1+r)^2} + \frac{1}{(1+r)^4}\frac{1- \(\frac{1}{1+r}\)^{d-3}}{1- \frac{1}{1+r}}-\frac{r}{(1+r)^3} \(  1-\(\frac r{1+r})^{d-3}\)\)\).
\cr}$$
Normalizing,
$$\eqalign{
\left\|\frac{1- r^{-1}X}{1+r^{-1}} \cdot \frac{Q}{(1+r)^d} \right\|  &  = \frac{r^2}{1+r} \(r^{-2}(1-r )+ r^{-1}\frac{1-r}{(1+r)^2}  +\frac{1}{r(1+r)^3}  - \frac{r}{(1+r)^4}\) + \Oh{\(\frac {r}{1+r}\)^d}\cr
& = \frac{1-r}{1+r}+ \frac {r(1-r)}{(1+r)^3}  + \frac{r}{(1+r)^4 } - \frac{r^3}{(1+r)^5} +  \Oh{\(\frac {r}{1+r}\)^d} \cr
& = \frac{1-r}{1+r} \(1 + \frac{1}{1+r} -\frac{1}{(1+r)^3}\)
+  \Oh{\(\frac {r}{1+r}\)^d} \cr
}$$

Now we note that if $P$ is a polynomial with degree $m$, and $Q$ is a polynomial of the form $q(x^M)$ where $M > m$, then $\| PQ\| = \| P \| \cdot \|Q \|$ (no mass cancellation can take place; in fact, if $x^a$ appears $PQ$ with nonzero coefficient, then there exists a unique pair $(c,d)$ \st $a= c+c'$, $c \in \Log P$, and $c' \in \Log Q$. 

We apply this with $P = \prod_{ j < k}P_j$ and $Q = p \prod_{k \leq j \leq d}$, where $p = (1+ rx^{2^k})/ (1+r)$ or $ (1+ r^{-1}x^{2^k})/ (1+r)^{-1} $. Therefore $s_{k,l} = (1-r)/(1+r)$ in the former case, and thus $\S((p_k), \NN (r)) = (1-r)/(1+r)$, whereas $\S((p_k\op), \NN(r))) = ((1-r)/(1+r) \cdot (1+ 1/(1+r) - 1/(1+r)^3) \neq \S((p_k), \NN(r))$. 
 \qed

This, together with the earlier results, yields $\NN(r) \iso \NN(r')$ implies $r= r'$ provided that $n(k) = n$ for all  $k$. 

   If  $n(k) \to \infty$ (the condition that merely $\sup n(k) = \infty$ appears to be much more complicated, and we do not deal with it), we run into a difficulty (although the mass-cancellation invariants can probably be used).  

 \Lem Proposition \twofiv. If $n(k) \to \infty$, then $\SS((w_k),\NN(r)) =1$ for all $r> 0$.

 \Pf The condition $n(k) \to \infty$ implies that for each $j$, the set  $\Set{l \in \N}{n(l) = j}$ is finite. It follows immediately that there exist infinitely many $l$ with the property that for all $k' > l$, we have $n(l+1) < n(k'+1)$. 

 It suffices, from the definition of $S_l$ to show that for every $k' > l $ (where $n(l+1) < n(k'+1)$ for all $k' > l$), that $S_{l+1,l} < S_{k',l}$; sufficient for this is, 
$$
 1 - \frac {4r \sin^2 \pi/(n(l+1))}{(1+r)^2} \leq  \prod_{t = 0}^{k'-l} \( 1 -  \frac{4r}{(1+r)^2} \sin^2 \frac  {\pi}{n(k'+1) n(k') \cdots n(l+1+t)}\). 
 $$ 
 To this end, we may assume that $l$ is so large that  for all $k > l$, we have $n(k+1) \geq 10$, and we observe that the left side is bounded below  by 
$$\eqalign{
 \prod_{t = 0}^{k'+l-1} \( 1 -  \frac{4r}{(1+r)^2} \(\frac{\pi}{n(k'+1) n(k') \cdots n(l+1+t)}\)^2\) &   \geq \(1 - \frac{4r \pi^2}{n(k'+1)^2 (1+r)^2}\) \times \cr & \qquad\(1 + \sum_{t=0}^{k'-l-1} \frac 1{(n(k') \cdot n(k'-t))^2} \)\cr
 &  \geq \( 1 -  \frac{4r \pi^2}{n(k'+1)^2 (1+r)^2}\) \times \cr
& \qquad \(1 + \frac{1}{(n(l+1)+1)^2 -1}\)\cr 
 }$$
 (The last line comes from $n(k'+1) > n(l+1)$.)

 Finally $1 - 4r \sin^2 (\pi/n(l+1)) /(1+r)^2 \leq 1- 4r \pi^2(1+ \eta) /n(k'+1)^2 (1+r)^2$ is equivalent to $ \sin^2 (\pi/n(l+1)) \geq (1 + \eta) \pi^2/n(k'+1)^2$. The latter is at least as large as $  \pi^2/((n(l+1)+1)^2-1)$. For $n(l+1)\geq 4$, we have $\sin^2 \pi/n(l+1) = (1-\cos 2\pi /n(l+1))/2  \geq \pi^2/n(l+1)^2 - 1\pi^4/3 n(l+1)^4$. So sufficient is that 
$$
 \frac{\pi^2}{n(l+1)^2} - \frac{\pi^4}{3n(l+1)^4} \geq \frac{\pi^2}{(n(l+1) +1)^2 -1}.
$$
 But this is a straightforward consequence of $6 n(l+1) \geq 10 > \pi^2$.\qed This finishes Example \twotwo. \qed 

Difficulties arise when we try to extend this to more general $h$. Suppose that  $h$ is a polynomial of degree $d > 1$, and as usual, $n(k) = n$ (constant). The behaviour of the product $S_l = \prod_{t=1}^{\infty} |h(r\exp (2\pi i/n^t))/h(r)|$ is   more complicated, when viewed as a function of $r$. Instead of having just one minimum value (as $r$ varies), it can have several critical points (up to $d$ of them). This can be somewhat compensated for.

To give an example, suppose that $h = (1 + x + 2x^2)/4$ and $n = 3$ (so we still have a non-interactive situation). Instead of taking $w_k = \exp ( 2\pi i/3^k)$, we may make another choice, $w_{k,2} = \exp ( 2\pi i/2\cdot 3^k)$. This yields another invariant, and the value will be not zero.   

For general polynomial $h$ and $n(k) = n$ (easier to deal with if $ n  > \deg h$), for $j=1,2,\dots d$,  we can use each of the sequences $(w_{k,j} =\exp  ( 2\pi i/j\cdot  n^k))_k$. This yields $d$ invariants, and a corresponding map $\R^{++} \to [0,1]^d$, given by $r \mapsto (\SS((w_{k,j}), \NN(r)))_{j=1}^d$. It is plausible that when $d > 1$ and $h$ is suitably nondegenerate, this is one to one---which would yield the nonisomorphism result.

\long\def\Rf[#1] #2, #3. #4\par%
{\vskip 2pt \itemitem{[#1]} #2, {\it #3,} #4\par\vskip2pt}

\SecT References

\Rf [BH] S Bezuglyi  and D Handelman, Measures on Cantor sets: the good, the ugly, the bad. Trans Amer Math Soc 366 (2014) 6247--6311.

\Rf [BD] R Bhatia and C Davis, A better bound on the variance. American Math Monthly  107 (2000)  353--357.

\Rf [CW] A Connes and EJ Woods, Hyperfinite von Neumann algebras and Poisson boundaries for time dependent random walks. Pacific J of Math  137 (1989)  225--243.

\Rf [EG] GA Elliott and T Giordano, Amenable actions of discrete groups. Erg Thy \& Dyn Sys 13 (1993) 289--318.

\Rf [GH] T Giordano and D Handelman, Matrix-valued random walks and variations on property AT. M\" unster  J  Math 1 (2008)  15--72. 

\Rf [H] D Handelman, Isomorphisms and nonisomorphisms of AT actions. J d'Analyse Math\'ematique 108 (2009)  293--396.

\Rf [H2] D Handelman, Nearly approximate transitivity (AT) for circulant matrices. Can J Math 71 (2019)  381--415.

\vskip 10pt

Mathematics Department, University of Ottawa, Ottawa ON  K1S 2H8, Canada; rochelle2\@sympatico.ca

\end